\documentclass[11pt,oneside]{article}
\usepackage{amssymb,amsmath}
\usepackage[english]{babel}
\usepackage{pstricks,pst-node}
\renewcommand{\le}{\leqslant}
\renewcommand{\ge}{\geqslant}

\newcommand{\chr}{\mathrm{char}}

\newcommand{\RR}{\mathbb{R}}
\newcommand{\ZZ}{\mathbb{Z}}
\newcommand{\QQ}{\mathbb{Q}}

\newcommand{\NN}{\mathbb{N}}
\newcommand{\FF}{\mathbb{F}}
\newcommand{\CC}{\mathbb{C}}

\newcommand{\VVV}{\mathcal{V}}

\newcommand{\vu}{\mathbf{u}}
\newcommand{\vv}{\mathbf{v}}

\newcommand{\ova}{\overline{a}}
\newcommand{\ovb}{\overline{\beta}}

\newtheorem{lemma}{Lemma}
\newtheorem{theorem}{Theorem}
\newtheorem{proposition}{Proposition}
\newtheorem{definition}{Definition}
\newtheorem{problem}{Problem}

\newtheorem{conjecture}[problem]{Conjecture}

\newtheorem{corollary}{Corollary}

\newtheorem{theorembhwy}{Theorem BHWY2}

\newtheorem{theorembhwyone}{Theorem BHWY1}

\newtheorem{theoreml}{Theorem L}

\newtheorem{theoremhone}{Theorem H1}

\newtheorem{theoremhtwo}{Theorem H2}

\newenvironment{proof}{\textsc{Proof. }}{\ \newline\hspace*{\fill}$\boxtimes$}

\newcommand{\bigk}{\mathop{\mathbf{K}}}

\begin{document}

\title{Continued fractions of certain Mahler functions}

\author{
 Dzmitry Badziahin
}

\maketitle

\begin{abstract}
We investigate the continued fraction expansion of the infinite
product $g(x) = x^{-1}\prod_{t=0}^\infty P(x^{-d^t})$ where the
polynomial $P(x)$ satisfies $P(0)=1$ and $\deg(P)<d$. We construct
relations between the partial quotients of $g(x)$ which can be used
to get recurrent formulae for them. We provide formulae for the
cases $d=2$ and $d=3$. As an application, we prove that for $P(x) =
1+ux$ where $u$ is an arbitrary rational number except 0 and 1, and
for any integer $b$ with $|b|>1$ such that $g(b)\neq0$ the
irrationality exponent of $g(b)$ equals two. In the case $d=3$ we
provide a partial analogue of the last result with several
collections of polynomials $P(x)$ giving the irrationality exponent
of $g(b)$ strictly bigger than two.
\end{abstract}

\section{Introduction}

Let $\FF$ be a field. Consider the set $\FF[[x^{-1}]]$ of Laurent
series together with the valuation which is defined as follows: for
$f(x) = \sum_{k=-d}^\infty c_kx^{-k} \in \FF[[x^{-1}]]$ its
valuation $||f(x)||$ is the biggest degree $d$ of $x$ having
non-zero coefficient $c_{-d}$. For example, for polynomials $f(x)$
the valuation $||f(x)||$ coincides with their degree. It is well
known that in this setting the notion of continued fraction is well
defined. In other words, every $f(x)\in \FF[[x^{-1}]]$ can be
written as
$$
f(x) = [a_0(x), a_1(x),a_2(x),\ldots],
$$
where the $a_i(x)$ are non-zero polynomials of degree at least 1,
$i\in\ZZ_{\ge 0}$. We provide some facts about the continued
fractions of Laurent series in Section~\ref{sec_contfrac} and refer
the reader to a nice survey~\cite{poorten_1998} for more details.

It appears that in the case $\FF=\QQ$ quite often the continued
fraction of $f(x)\in \QQ[[x^{-1}]]$ can give us the information
about the approximational properties of real numbers $f(b)$ for
integer values $b$ inside the radius of convergence of $f(x)$. One
of the most important such properties is an irrationality exponent.
It indicates how well a given irrational number $\xi$ is
approximated by rationals and is denoted by $\mu(\xi)$. More
precisely, it is the supremum of real numbers $\mu$ such that the
inequality
$$
\left|\xi - \frac{p}{q}\right|<q^{-\mu}
$$
has infinitely many rational solutions $p/q$. This is one of the
most important approximational properties of real numbers indicating
how well $\xi$ is approximated by rationals. Note that by the
classical Dirichlet approximation theorem we always have
$\mu(\xi)\ge 2$.

%
%

Let $f(x)$ be an infinite product defined as follows
\begin{equation}\label{main_eq}
f(x) = \prod_{t=0}^\infty P(x^{-d^t})
\end{equation}
where $P\in \FF[x]$ is a polynomial and $d\ge 2$ is a positive
integer. In order that $f(x)$ is correctly defined as a Laurent
series, we need an additional condition $P(0)=1$. An easy check
shows that functions $f(x)$ fall into the set of {\em Mahler
functions} which we define as follows: $M(x)\in\FF[[x^{-1}]]$ is a
Mahler function if it satisfies the equation of the form
$$
\sum_{i=0}^n P_i(x)M(x^{d^i})=0
$$
for some integers $n\ge 1, d\ge 2$, and polynomials $P_0(x),\ldots, P_n(x)\in
\FF[x]$ with $P_0(x)P_n(x)\neq 0$\footnote{In the literature the notion of
Mahler function is often given to $F(x)\in \QQ[[x]]$ and not to Laurent
series $M(x)$ as in our case. However one can easily convert one notion into
another by considering $M(x) = F(x^{-1})$.}. For any integer $b$ within the
radius of convergence of $M$ the value $M(b)$ is called {\em Mahler number}.

The question about computing or at least estimating the
irrationality exponent of Mahler numbers is currently in the focus
of the Diophantine approximation. It was triggered by the work of
Bugeaud~\cite{bugeaud_2011} where he showed that for $b\ge 2$ the
irrational exponent of the Thue-Morse numbers $f_{TM}(b)$ is equal
to 2. Here $f_{TM}(x)$ is the most classical example of Mahler
functions and can be defined as follows:
$$
f_{TM}(x) := \prod_{t=0}^\infty (1 - x^{-2^t}).
$$
One of the key ingredients of that paper is the result
from~\cite{alouche_peyriere_wen_wen_1998} about non-vanishing of
Hankel determinants of $f_{TM}(x)$ (they will be properly defined
and discussed in Section~\ref{sec_hankel}). Later this approach was
further developed and generalised to cover many other Mahler
functions, see for example~\cite{coons_2013, guo_wu_wen_2014,
wu_wen_2014}. Finally, Bugeaud, Han, Wen and Yao
\cite{bugeaud_han_wen_yao_2015} provided quite a general result
where the estimates for $\mu(f(b))$ are given depending on the
distribution of non-vanishing Hankel determinants of $f(x)$ (see
Theorem BHWY2 in Section~\ref{sec_hankel}). The problem with this
theorem is that it is usually quite problematic to compute the
Hankel determinants of $f(x)$ or even to check which of them is
equal to zero. In \cite{han_2015a, han_2015b,
bugeaud_han_wen_yao_2015} the authors used the reduction of $f(x)$
modulo a prime number $p$ to provide local conditions on $f(x)$
which ensure that $\mu(f(b)) = 2$. We present just one example of
such results, which appears as Theorem~2.5 in
\cite{bugeaud_han_wen_yao_2015}.

\begin{theorembhwyone}
Let $f(x)\in\ZZ[[x^{-1}]]$ be a power series defined by
$$
f(x) = \prod_{t=0}^\infty \frac{C(x^{-3^t})}{D(x^{-3^t})},
$$
with $D(x),C(x)\in\ZZ[x]$ such that $C(0)=D(0)=1$. Let $b\ge 2$ be
an integer such that $C(b^{-3^m})D(b^{-3^m})\neq 0$ for all integer
$m\ge 0$. If $f(x) \pmod3$ is not a rational function then the
irrationality exponent of $f(b)$ is equal to 2.
\end{theorembhwyone}

Theorem~BHWY1 as well as other known results of this kind provide
infinite collections of Mahler functions $f(x)$ such that their
values $f(b)$ have irrationality exponent equal to 2. However,
firstly, many series $f(x)$ are not covered by the reduction modulo
$p$ approach. Secondly, it can not detect the cases when the
irrationality exponent of $f(b)$ is strictly bigger than two.

In Section~\ref{sec_hankel} we show that values of the Hankel
determinants of $f(x)$ can be derived from the continued fraction of
$f(x)$. Therefore in view of Theorem~BHWY2, understanding the
continued fraction gives us a powerful tool to estimate the
irrationality exponent of $f(b)$.

The question of computing the continued fraction of certain Mahler
functions (to the best of authors knowledge) goes back to 1991, when
Allouche, Mend\`es France and van der
Poorten~\cite{allouche_france_poorten_1991} showed that all partial
quotients of the infinite product
$$
f_3(x) := \prod_{t=0}^\infty (1-x^{-3^t})
$$
are linear. In~\cite{poorten_shallit_1992} the author computed the
continued fraction of the solution $f_M(x)$ of the equation
$$
f_M(x^2) = xf_m(x) - x^3.
$$
Some other papers on this topic are~\cite{poorten_1993,
poorten_1998}. In particular, in the second of these papers, van der
Poorten noted that the continued fraction of the Thue-Morse series
$$
f_2(x):=\prod_{t=0}^{\infty} (1-x^{-2^t})
$$
has a regular structure. In~\cite{badziahin_zorin_2014} the precise
formula of the continued fraction of $x^{-1}f_2(x)$ was given. As a
consequence of that the authors showed that the Thue-Morse constant
$f_2(2)$ is not badly approximable.
Later~\cite{badziahin_zorin_2015} the authors extended their result
to the series
$$
f_d(x):=\prod_{t=0}^{\infty} (1-x^{-d^t})
$$
for any $d\ge 2$. In particular, they show that $f_d(x)$ is badly
approximable only for $d=2$ and $d=3$ (definition will be given in
Section~\ref{sec_contfrac}) and provide the formula for the
continued fraction of $x^{-2}f_3(x)$.

\subsection{Main results}

In this paper we consider functions $g(x) = x^{-1}f(x)$, where
$f(x)$ is given by an infinite product~\eqref{main_eq}. The
essential restriction we have to impose on them is $d>||P(x)||$
because that allows us, given a convergent of $g(x)$, to produce an
infinite chain of other convergents of $g(x)$. Under these
restrictions we encode each function $g(x)$ by a vector
$\vu=(u_1,\ldots,u_{d-1})\in\FF^{d-1}$ in the following way:
\begin{equation}\label{main_ueq}
g_\vu(x):= x^{-1}\prod_{t=0}^\infty
(1+u_1x^{-d^t}+u_2x^{-2d^t}+\ldots+u_{d-1}x^{-(d-1)d^t}).
\end{equation}
The notation $f_\vu(x)$ is defined in the same way.

We managed to find the relations between the partial quotients of
the continued fraction of $g_\vu(x)$ which can provide the recurrent
formulae for them (see Propositions~\ref{prop1} and~\ref{prop2} in
Section~\ref{sec_conds}). We then explicitly write down these
recurrent formulae in the case $d=2$ and $d=3$:

\begin{theorem}\label{th1}
Let $u\in\FF$. If $g_u(x)$ is badly approximable then its continued
fraction is
$$
g_u(x) = \bigk_{i=1}^\infty \frac{\beta_i}{x+\alpha_i}
$$
where the coefficients $\alpha_i$ and $\beta_i$ are computed by the
formula
\begin{equation}\label{recur_d2}
\begin{array}{l}
\alpha_{2k+1} = -u,\; \alpha_{2k+2}=u;\\
\displaystyle\beta_1 = 1,\; \beta_2 = u^2-u,\; \beta_{2k+3} =
-\frac{\beta_{k+2}}{\beta_{2k+2}},\; \beta_{2k+4} = \alpha_{k+2}+u^2
- \beta_{2k+3}
\end{array}
\end{equation}
for any $k\in\ZZ_{\ge 0}$.
\end{theorem}

\begin{theorem}\label{th2}
Let $\vu=(u,v)\in\FF^2$. If $g_\vu(x)$ is badly approximable then
its continued fraction is
$$
g_u(x) = \bigk_{i=1}^\infty \frac{\beta_i}{x+\alpha_i}
$$
where the coefficients $\alpha_i$ and $\beta_i$ are computed by the
recurrent formula
\begin{equation}\label{init_d3}
\begin{array}{lll}
\alpha_1 = -u,&\displaystyle \alpha_2 =
\frac{u(2v-1-u^2)}{v-u^2},&\displaystyle\alpha_3 =
\frac{-u(v-1)}{v-u^2};\\[2ex]
\beta_1 = 1,& \beta_2 = u^2-v,&\displaystyle\beta_3 =
\frac{u^2+u^4+v^3 - 3u^2v}{(v-u^2)^2}.
\end{array}
\end{equation}
and
\begin{equation}\label{recur_d3}
\begin{array}{l}
\displaystyle \alpha_{3k+4} = -u,\quad \beta_{3k+4} =
\frac{\beta_{k+2}}{\beta_{3k+3}\beta_{3k+2}};\\[1ex]
\displaystyle \beta_{3k+5} = u^2 - v - \beta_{3k+4},\quad
\alpha_{3k+5} = u -
\frac{\alpha_{k+2}+uv-\alpha_{3k+2}\beta_{3k+4}}{\beta_{3k+5}}\\[1ex]
\alpha_{3k+6} = u-\alpha_{3k+5},\quad \beta_{3k+6} = v -
\alpha_{3k+5}\alpha_{3k+6}.
\end{array}
\end{equation}
for any $k\in\ZZ_{\ge 0}$.
\end{theorem}

As one can notice, the complexity of the recurrent formulae grows
rapidly with $d$. Therefore a computer assistance may be needed to
provide analogues of Theorems~\ref{th1} and~\ref{th2} for larger
values of $d$. These two theorems are proven in
Section~\ref{sec_conds}.

We show that the formulae in Theorems~\ref{th1} and~\ref{th2} can be
used to check whether $g_\vu(x)$ is badly approximable or not:

\begin{theorem}\label{th7}
Let $\FF\subset \CC$. The function $g_u(x)$ (respectively,
$g_\vu(x)$) is badly approximable if and only if none of the
parameters $\beta_n$, computed by formulae~\eqref{recur_d2}
(respectively~\eqref{init_d3} and~\eqref{recur_d3}) vanish.
Moreover, if $\beta_1,\beta_2,\ldots, \beta_n\neq 0$ then the first
$n$ partial quotients of $g_u(x)$ (respectively $g_\vu(x)$) are
linear. And if in addition $\beta_{n+1}=0$ then the $(n+1)$th
partial quotient of $g_u(x)$ is not linear.
\end{theorem}

Theorem~\ref{th7} is probably true for any field $\FF$, however, as
we are mostly interested in $\FF=\QQ$, we proved it only for
subfields of complex numbers and did not make a big effort to
generalise the result to an arbitrary field. We prove this result in
the beginning of Subsection~\ref{subsec_bad23}.

For the remaining results we assume that $\FF=\QQ$. Equipped with
the continued fraction of $g_\vu(x)\in\QQ[[x^{-1}]]$ we can compute
or at least estimate irrationality exponents of the values
$g_\vu(b)$. The first result we want to provide here is as follows:

\begin{theorem}\label{th8}
Let $b\ge 2$ be integer such that $g_\vu(b)\neq 0$. Then
$\mu(g_\vu(b))=2$ if and only if $g_\vu(x)$ is badly approximable.
\end{theorem}

The ``if'' part of this theorem is covered in
Section~\ref{sec_hankel}. As we will see, it is essentially an
implication of Theorem~BHWY2 from~\cite{bugeaud_han_wen_yao_2015}.
The ``only if'' part is considered in Section~\ref{sec_prelim}.

In the case when $g_\vu(x)$ is not badly approximable we provide a
non-trivial lower bound for the irrationality exponent of
$\mu(g_\vu(b))$. This result is proven in Section~\ref{sec_irr}.

\begin{theorem}\label{th5}
Let $b\ge 2$ be integer. If $g_\vu(x)$ is not badly approximable
 and $g_\vu(b)\neq 0$ then the irrationality exponent of $g_\vu(b)$ satisfies
$$
\mu(g_\vu(b))\ge 2+\frac{c-1}{n_0},
$$
where $n_0$ is the smallest positive value such that the $(n_0+1)$-th
convergent $a_{n_0+1}(x)$ is not linear and $c = ||a_{n_0+1}(x)||$.
\end{theorem}

We finish the paper by applying the results from above to compute
(or estimate) the irrationality exponents of $g_{\vu}(b)$ for all
integer values $b\ge 2$ and as many vectors $\vu$ as possible. We
manage to completely cover the case $d=2$ and $u\in\QQ$:

\begin{theorem}\label{th3}
The series $g_u(x)$ is badly approximable for any $u\in\QQ$ except
$u=1$ and $u=0$ for which $g_u(x)$ becomes a rational function. In
particular, if $u\in \QQ\backslash\{0,1\}, b\in\ZZ$, $|b|>1$ and
$b^{2^t} + u\neq 0$ for any $t\in\ZZ_{\ge 0}$ then the real number
$g_u(b)$ has irrationality exponent two.
\end{theorem}

However, due to the complexity of the formulae~\eqref{init_d3}
and~\eqref{recur_d3} we covered many but not all values of
$\vu\in\QQ^2$ for $d=3$.

\begin{theorem}\label{th4}
The series $g_{(u,0)}(x)$ as well as $g_{(0,v)}(x)$ is badly
approximable for all $u,v\in\QQ$ except $u=0$ or $v=0$ respectively.
In particular, if $u\in \QQ\backslash\{0\}, b\in\ZZ, |b|>1$ then the
irrationality exponent of $g_{(u,0)}(b)$ and of $g_{(0,v)}(x)$ is
two as soon as $b^{3^t}+u\neq 0$ and $b^{2\cdot 3^t}+v\neq 0$ for
any $t\in \ZZ_{\ge 0}$.
\end{theorem}

\begin{theorem}\label{th6}
Let $\vu=(u,v)\in\RR^2$ satisfy the following conditions:
\begin{enumerate}
\item[(C1)] $u^2\ge 6$;
\item[(C2)] $v\ge \max\{3u^2-1, 2u^2+8\}$.
\end{enumerate}
Then the series $g_\vu(x)$ is badly approximable. In particular, if
$\vu\in\QQ^2$, $b\in\ZZ$ satisfies $|b|>1$ and $b^{2\cdot
3^t}+ub^{3^t}+v\neq 0$ for any $t\in\ZZ_{\ge 0}$, the irrationality
exponent of $g_\vu(x)$ is two.
\end{theorem}

In the proof of Theorem~\ref{th6} we sometimes make quite rough
estimates, therefore with no doubts the conditions~(C1) and~(C2) can
be made weaker. By this theorem we want to demonstrate that the
knowledge of the continued fraction of $g_\vu(x)$ can produce global
conditions on $\vu$ for the series to be badly approximable, on top
of the local conditions, as in Theorem~BHWY1.

Finally, by investigating the equations $\beta_n=0$ for small values
of $n$ we get several series of vectors $\vu\in\QQ^2$ such that
$g_\vu(x)$ is not badly approximable and therefore non-trivial lower
bounds for $g_\vu(b)$ apply:

\begin{theorem}\label{th9}
The functions $g_\vu(x)$ are not badly approximable for the
following vectors $\vu\in\QQ^2$:
\begin{enumerate}
\item $\vu = (\pm u,u^2)$. Then for any $u\in\QQ$ and any $b\in\ZZ$ with
    $|b|>1$ we have $\mu(g_\vu(b))\ge 3$, as soon as
    $g_\vu(x)\not\in\QQ(x)$ and $g_\vu(b)\neq 0$.
\item $\vu = (\pm s^3, -s^2(s^2+1))$. Then for any $s\in\QQ$ and any
    $b\in\ZZ$ under the same conditions as before we have $\mu(g_\vu(b))\ge
    3$.
\item $\vu = (2,1)$. Then for any $b\in\ZZ$ with $|b|>1$ we have
    $\mu(g_\vu(b))\ge 12/5$.
\end{enumerate}
\end{theorem}

We wrote a computer program which computed the first 30 partial
quotients for all integer values $u,v$ in the range $|u|,|v|<1000$
and it did not find any other values $\vu= (u,v)$ than those
mentioned above, for which $g_\vu(x)$ is not badly approximable.
This suggests that the following statement may take place:

\begin{conjecture}
The only values $\vu\in\ZZ^2$ such that $g_\vu(x)$ is not badly
approximable are as follows:
$$
(\pm u,u^2),\; (\pm s^3,-s^2(s^2+1))\;\;\mbox{and}\;\; (\pm2,1),
$$
where $s\in\ZZ$.
\end{conjecture}

By observing that $g_{(2,1)}(x) = (g_{1,0}(x))^2$ we get a notable
corollary from Theorems~\ref{th4} and~\ref{th9}: it provides a
family of Mahler numbers $\xi$ such that $\mu(\xi)=2$, but
$\mu(\xi^2)\ge 7/3$. Namely, one can take $\xi = g_{1,0}(b)$ for any
integer $b$ with $|b|>1$.

\noindent{\bf Remark.} For any integer $b$, Mahler numbers
$g_\vu(b)$ and $f_\vu(b)$ are rationally dependent. Therefore they
share the same irrationality exponent and Theorems~\ref{th8} --
\ref{th9} also provide the information about the irrationality
exponents of perhaps ``nicer'' Mahler numbers $f_\vu(x)$. Also, as
explained in in Section~\ref{sec_contfrac},
Theorems~\ref{th7},~\ref{th3} -- \ref{th9} give us an insight
whether the function $f_\vu(x)$ is badly approximable or not.
However its continued fraction definitely differs from what is
provided in the first two theorems.

\section{Continued fractions and continuants}\label{sec_contfrac}

Continued fractions of Laurent series share many properties of the
classical continued fractions in real numbers. For example, it is
known that, as for the standard case, the convergents $p_n(x)/q_n(x)
= [a_0(x);a_1(x),\ldots, a_n(x)]$ of $f(x)$ are the best rational
approximants of $f(x)$. Furthermore, we have an even stronger
version of Legendre's theorem:
\begin{theoreml}\label{th_legendre}
Let $f(x)\in \FF[[x^{-1}]]$. Then $p(x)/q(x)\in \FF(x)$ in a reduced
form is a convergent of $f(x)$ if and only if
$$
\left|\left| f(x) - \frac{p(x)}{q(x)}\right|\right| < -2||q(x)||.
$$
\end{theoreml}
The proof of this and other unproven facts from this section can be
found in~\cite{poorten_1998}.

As we already mentioned, every series $f(x)\in\FF[[x^{-1}]]$ allows
an expansion into a continued fraction. We will use the following
notation:
$$
f(x) := [a_0(x);a_1(x), a_2(x),\ldots] = a_0(x) + \bigk_{n=1}^\infty
\frac{1}{a_n(x)},
$$
where $a_i(z)\in \FF[z], i\in\NN$. The convergents $p_n(x)/q_n(x)$
of $f(x)$ can be computed by the following formulae
\begin{equation}\label{conv_def}
\begin{array}{l}
p_{-1}(x) = 1,\; p_0(x) = a_0(x),\; p_{n+1}(x) = a_{n+1}(x)p_n(x) +
p_{n-1}(x);\\[1ex] q_{-1}(x) = 0,\; q_0(x) = 1,\; q_{n+1}(x) =
a_{n+1}(x)q_n(x) + q_{n-1}(x).
\end{array}
\end{equation}
However, unlike the classical setup of real numbers, where the
numerators and the denominators $p_n$ and $q_n$ are defined
uniquely, $p_n(x)$ and $q_n(x)$ are only unique up to multiplication
by a non-zero constant. We can make them unique by putting an
additional condition, that $q_n(x)$ must be monic. It is not
difficult to see that~\eqref{conv_def} do not usually give monic
polynomials. However we can adjust these formulae a bit to meet the
required condition:
\begin{equation}\label{mconv_def}
\begin{array}{ll}
\hat{p}_{-1}(x) = 1,\;\; \hat{p}_0(x) = a_0(x),& \hat{p}_{n+1}(x) =
\hat{a}_{n+1}(x)\hat{p}_n(x) + \beta_{n+1}
\hat{p}_{n-1}(x);\\[1ex] \hat{q}_{-1}(x) = 0,\;\; \hat{q}_0(x) = 1,& \hat{q}_{n+1}(x) =
\hat{a}_{n+1}(x)\hat{q}_n(x) + \beta_{n+1}\hat{q}_{n-1}(x).
\end{array}
\end{equation}
where we define, with $\rho_n$ denoting the leading coefficient of
$q_n(x)$:
$$
\hat{a}_{n+1}(x) := \frac{a_{n+1}(x)\cdot
\rho_n}{\rho_{n+1}}\quad\mbox{and}\quad \beta_0 = \beta_1=1,\;
\beta_{n+1} = \frac{\rho_{n-1}}{\rho_{n+1}}.
$$
One can easily check from~\eqref{mconv_def} that $\hat{a}_n(x)$ are
always monic. The formula for $\beta_n$ suggests that $\beta_m\neq
0$.

It is not difficult to check that from the sequence of monic
polynomials $\hat{a}_n(x)$ together with the sequence of non-zero
elements $\beta_n$ one can uniquely restore the initial continued
fraction $[a_0(x),a_1(x),\ldots]$. Indeed, we have $a_n(x) =
\rho_n\hat{a}_n(x)$ and $\rho_n$ can be derived from the formula
$\beta_{n+1} = \rho_{n-1}/\rho_{n+1}$ and initial values
$\rho_0=\rho_1=1$. In other words, any Laurent series has a modified
continued fraction of the form
\begin{equation}\label{cfrac}
f(x) = \hat{a}_0(x) +
\frac{\beta_1}{\hat{a}_1(x)+\frac{\beta_2}{\hat{a}_2(x)+\cdots}} =:
\hat{a}_0(x) + \bigk_{n=1}^\infty \frac{\beta_n}{\hat{a}_n(x)},
\end{equation}
where $\hat{a}_n(x)\in\FF[x]$ are monic and $\beta_n\in\FF$ are
non-zero. And vice versa: any sequence of monic $\hat{a}_n(x)$ and
non-zero values $\beta_n$ defines a modified continued fraction for
some $f(x)$.

In the paper we will use the modified formulae~\eqref{mconv_def} for
computing convergents and also for convenience we will not write
hats above the variables $p_n,q_n$ and $a_n$.

For our function $g(x) = x^{-1}f(x)$ where $f(x)$ is defined
by~\eqref{main_eq}, we have $x^{-1} + u_1x^{-2}+\ldots$, where $u_1$
is the coefficient coming from $P(x) = 1+u_1x+\ldots$. Therefore the
first values of the convergents for $g(x)$ are computed as follows:
$$
\begin{array}{l}
p_{-1}(x) = 1,\;\;p_0(x) = 0,\;\; p_1(x) = 1;\\
q_{-1}(x) = 0,\;\;q_0(x) = 1,\;\; q_1(x) = x-u_1 = a_1(x).
\end{array}
$$
Notice that in this case formulae~\eqref{mconv_def} give both monic
$p_n(x)$ and $q_n(x)$.

Equations \eqref{mconv_def} can be written in terms of the
generalised continuants as follows. Given two sequences
$\overline{a} = (a_n(x))_{n\in\NN}$ and $\overline{\beta} =
(\beta_n)_{n\in\NN}$ and $k\le l$ we write
$$
\begin{array}{l}
\overline{a}_{k,l}:= (a_k(x),a_{k+1}(x),\ldots,a_l(x)),\quad
\overline{a}_{k} := \overline{a}_{1,k};\\[1ex]
\overline{\beta}_{k,l} := (\beta_k,\beta_{k+1},\ldots,\beta_l),\quad
\overline{\beta}_k:= \overline{\beta}_{1,k}.
\end{array}
$$
Now generalised continuants $K_n^0(\ova_n, \ovb_n)\in\FF[x]$ and
$K_n^1(\overline{a}_n, \overline{\beta}_n)\in\FF[x]$ are defined as
follows: $K^0_{-1} () = 0,\;K^1_{-1}()=1;\;K^0_0() = 1;\; K^1_0()=0$
and for $n\in\NN$ both $K^0_n(\ova_n,\ovb_n)$ and
$K^1_n(\ova_n,\ovb_n)$ satisfy the same recurrent equation:
\begin{equation}\label{cont_def}
K_n(\ova_n,\ovb_n) =
a_n(x)K_{n-1}(\ova_{n-1},\ovb_{n-1})+\beta_nK_{n-2}(\ova_{n-2},\ovb_{n-2}).
\end{equation}

One can easily check that the $K_n^0(\ova_n,\ovb_n)$ are always
monic while the leading coefficient of $K_n^1(\ova_n,\ovb_n)$ equals
$\beta_1$ for any $n\ge 1$. We can also check that the degrees of
the continuants satisfy
\begin{equation}\label{cont_deg}
||K_n^0(\ova_n,\ovb_n)|| = \sum_{k=1}^n ||a_k(x)||,\quad
||K_n^1(\ova_n,\ovb_n)|| = \sum_{k=2}^n ||a_k(x)||.
\end{equation}

The enumerator and the denominator of the $n$'th convergent of
$g(x)$ can be written as $p_n(x) = K^1_n(\ova_n,\ovb_n)$ and $q_n(x)
= K^0_n(\ova_n,\ovb_n)$. Moreover, these polynomials are linked
together by the following relations:

\begin{lemma}\label{lem7}
For $n,m\in\ZZ_{\ge 0}$ we have
\begin{equation}\label{lem7_eq}
p_{n+m}(x) = K_m^0(\ova_{n+1,n+m},\ovb_{n+1,n+m})p_n(x) +
K_m^1(\ova_{n+1,n+m},\ovb_{n+1,n+m})p_{n-1}(x).
\end{equation}
The same relation is true for the polynomials $q_n(x)$ too.
\end{lemma}
Formula~\eqref{lem7_eq} can be checked by applying~\eqref{mconv_def}
and by using induction on $m$.

We will need to quantify the inequality from Theorem~L.

\begin{definition}
Let $p(x)/q(x) \in \FF(x)$ be a rational function and $u(x)$ be a
Laurent series. We say that an integer $c>0$ is {\em the rate of
approximation of $u(x)$ by $p(x)/q(x)$} if
$$
||u(x) - p(x)/q(x)|| = -2||q(x)||-c
$$
\end{definition}

It is known (see \cite[displayed equation before Proposition
1]{poorten_1998}) that the convergent $p_n(x)/q_n(x)$ approximates
$u(x)$ with the rate $||a_{n+1}(x)||$.

\begin{definition}
We say that $f(x)\in \FF((x^{-1}))$ is {\em badly approximable} if
each valuation (i.e. degree) of its partial quotients is bounded
from above by an absolute constant. Otherwise we say that $f(x)$ is
well approximable.
\end{definition}

An equivalent formulation of this definition is: $f(x)$ is badly
approximable if for any $n\in\NN$ the rate of approximation of $f$
by its $n$th convergent is bounded from above by an absolute
constant.

We end this section with a lemma which shows that the following two
statements are equivalent: $f_\vu(x)$ is badly approximable and
$g_\vu(x)$ is badly approximable. Its proof can be found
in~\cite[Proposition 1]{badziahin_zorin_2015}.

\begin{lemma}
Let $f(x)\in\QQ[[x^{-1}]]$, $a(x), b(x)\in\QQ[x]\backslash {0}$.
Then $f(x)$ is badly approximable if and only if $g(x) =
\frac{a(x)}{b(x)}f(x)$ is badly approximable.
\end{lemma}

\section{Relation with Hankel continued fractions}\label{sec_hankel}

As we mentioned in the Introduction, the more popular approach to
compute irrationality exponents of Mahler numbers uses Hankel
determinants and Hankel continued fractions rather than the
classical ones. For example, they can be found in the works of
Han~\cite{han_2015a,han_2015b}. For the power series
$f(x)\in\FF[[x]]$ they are defined as follows
\begin{equation}\label{hank}
f(x) =
\frac{v_0x^{k_0}}{1+u_1(x)x-\frac{v_1x^{k_0+k_1+2}}{1+u_2(x)x-\frac{v_2x^{k_1+k_2+2}}{1+u_3(x)x-\ldots}}}
\end{equation}
where $v_i\neq 0$ are constants, $k_i$ are nonnegative integers and
$u_i(x)$ are polynomials of degree $\deg(u_i)\le k_{i-1}$. For
convenience we will use the following shorter notation instead
of~\eqref{hank}:
$$
f(x) = \bigk_{i=1}^\infty \frac{v^*_{i-1}
x^{k_{i-1}+k_{i-2}+2}}{1+u_{i}(x)x}.
$$
where $v^*_0 = v_0$ and $v^*_i = -v_i$ for $i\ge 1$.

In particular, the following result was established
in~\cite{han_2015b}:

\begin{theoremhone}
Each Hankel continued fraction defines a power series and
conversely, for each power series $f(x)$ there exists unique Hankel
continued fraction of $f(x)$.
\end{theoremhone}

If we consider the Laurent series $x^{-1}f(x^{-1})$ then the Hankel
continued fraction transforms to the standard continued fraction in
the space of Laurent series. Indeed, one can easily check that
$$
x^{-1}f(x^{-1}) = \bigk_{i=1}^\infty
\frac{v_{i-1}}{x^{k_{i-1}+1}+x^{k_{i-1}}u_i(x^{-1})}.
$$
Notice that if we set $a_i(x) =
x^{k_{i-1}+1}+x^{k_{i-1}}u_i(x^{-1})$ and $\beta_i = v_{i-1}$ then
we get the same notation as in~\eqref{cfrac}. Since $a_i(x)$ is
surely a polynomial, this gives a one-to-one correspondence between
Hankel continued fractions for $f(x)$ and standard continued
fractions for Laurent series of $x^{-1}f(x^{-1})$. In particular,
this observation together with the standard fact that continued
fractions for $x^{-1}f(x)$ are uniquely defined, gives another proof
of Theorem~H1.

By applying~\eqref{conv_def} one can get that the degree of the
denominator $q_n(x)$ of the $n$th convergent of $x^{-1}f(x^{-1})$
can be computed as follows:
\begin{equation}\label{def_degq}
||q_n(x)|| = s_n = k_0+k_1+\ldots+k_{n-1}+n.
\end{equation}

Han used Hankel continued fractions of $f(x)$ to extract some
information about the Hankel determinants of $f(x)$ which prove to
be a powerful tool for computing the irrationality exponent of
numbers $f(b^{-1})$ where $b$ is a positive integer such that
$b^{-1}$ is inside the radius of convergence of $f(x)$. For $f(x) =
\sum_{n=0}^\infty c_nx^n$ the Hankel determinants are defined as
follows
$$
H_n(f) := \left|\begin{array}{cccc} c_0&c_{1}&\cdots
&c_{n-1}\\
c_{1}&c_{2}&\cdots&c_{n}\\
\vdots&\vdots&\ddots&\vdots\\
c_{n-1}&c_{n}&\cdots&c_{2n-2}
\end{array}\right|.
$$

The following result from~\cite{han_2015b} gives the relation
between the Hankel continued fraction and the Hankel determinants:

\begin{theoremhtwo}\label{th_hankel}
Let $f(x)$ be a power series such that its Hankel continued fraction
is given by~\eqref{hank}. Then, for all integers $i\ge 0$ all
non-vanishing Hankel determinants are given by
\begin{equation}\label{hankdet}
H_{s_i}(f) = (-1)^\epsilon v_0^{s_i}v_1^{s_i-s_1}\cdots
v_{i-1}^{s_i-s_{i-1}},
\end{equation}
where $\epsilon=\sum_{j=0}^{i-1} k_j(k_j+1)/2$ and $s_i =
k_0+k_1+\cdots+k_{i-1}+i$.
\end{theoremhtwo}

By combining this theorem with~\eqref{def_degq} we get
\begin{corollary}\label{cor2}
The $n$'th Hankel determinant of $f(x)$ does not vanish if and only
if there exists a convergent $p(x)/q(x)$ of $x^{-1}f(x^{-1})$ such
that $p(x)$ and $q(x)$ are coprime and $||q(x)||=n$.
\end{corollary}

Another straightforward corollary of the Theorem, applied to
continued fractions of $x^{-1}f(x^{-1})$ is as follows:
\begin{corollary}\label{cor3}
If the continued fraction of $x^{-1}f(x^{-1})$ is badly approximable
then there exists an increasing sequence $(s_i)_{i\ge 0}$ of
positive integers such that $H_{s_i}(f)\neq 0$ for all $i\in
\ZZ_{\ge 0}$ and $s_{i+1}-s_i$ is bounded from above by an absolute
constant dependent only on $f(x)$.
\end{corollary}

Guo, Wu and Wen~\cite{guo_wu_wen_2014} discovered a relation between
the sequence $s_i$ and the irrationality exponent of $f(b^{-1})$ for
Mahler functions $f(x)$. Their result was significantly improved and
corrected by Bugeaud, Han, Wen and
Yao~\cite{bugeaud_han_wen_yao_2015}.

\begin{theorembhwy}\label{th_bhwy}
Let $d\ge 2$ be an integer and $f(x) = \sum_{n=0}^\infty c_nx^n$
converge inside the unit disk. Suppose that there exist integer
polynomials $A(x), B(x), C(x), D(x)$ with $B(0)D(0)\neq 0$ such that
\begin{equation}\label{mahl_eq}
f(x) = \frac{A(x)}{B(x)} + \frac{C(x)}{D(x)}f(x^d).
\end{equation}
Let $b\ge 2$ be an integer such that
$B(b^{-d^n})C(b^{-d^n})D(b^{-d^n})\neq 0$ for all $n\in\ZZ_{\ge 0}$.
If there exists an increasing sequence $(s_i)_{i\ge 0}$ of positive
integers such that $H_{s_i}(f)\neq 0$ for all $i\in \ZZ_{\ge 0}$ and
$\limsup_{i\to \infty}\frac{s_{i+1}}{s_i}=\rho$, then $f(1/b)$ is
transcendental and
$$
\mu(f(1/b))\le (1+\rho)\min\{\rho^2,d\}.
$$
\end{theorembhwy}

We apply this theorem to the series
$$
\hat{f}_\vu(x):= f_\vu(x^{-1}) = \prod_{t=0}^\infty P(x^{d^t}).
$$
Note that it satisfies the equation $\hat{f}_\vu(x) =
\frac{1}{P(x)}\hat{f}_\vu(x^d)$ which is of the
form~\eqref{mahl_eq}. The condition
$B(b^{-d^n})C(b^{-d^n})D(b^{-d^n})\neq 0$ in this setting basically
means that $g_\vu(x)\neq 0$. Finally, an application of
Corollary~\ref{cor3} asserts that if $g_\vu(x)$ is badly
approximable and $g_\vu(b)\neq 0$ then $\mu(g_\vu(b))=2$. This
finishes the proof of ``if'' part of Theorem~\ref{th8}.

\noindent{\bf Remark.} The condition $s_{i+1}-s_i\le C$ is much
stronger than $\limsup_{i\to \infty} \frac{s_{i+1}}{s_i} = 1$, thus
the natural question arises: can we say anything better about the
approximational properties of $g_\vu(b)$ in the case $g_\vu(x)$ is
badly approximable? For example, can we show that
\begin{equation}\label{log_est}
\left|f\left(\frac{1}{b}\right) - \frac{p}{q}\right| \le
\frac{1}{q^2\cdot \delta(q)}
\end{equation}
for some $\delta(q)$ which grows slower than any power function
$q^\epsilon$? It appears that the proof
in~\cite{bugeaud_han_wen_yao_2015} can not be easily improved to
give us anything like~\eqref{log_est}.

\section{Information about $g(x)$}\label{sec_prelim}

Recall that we are focused on the following function written as a
Laurent series:
$$
g_\vu(x) = x^{-1}f_\vu(x) = x^{-1}\prod_{t=0}^{\infty} P(x^{-d^t}),
$$
where $P(x)\in \FF[x]$, $\deg(P)\le \dim(\vu) = d-1$ and $P(0) = 1$.
By substituting $x^d$ into the formula instead of $x$ we get the
functional relation
\begin{equation}\label{eq_g}
P^*(x) g_\vu(x^d) = g_\vu(x),
\end{equation}
where
\begin{equation}\label{pstar_def}
P^*(x) = x^{d-1}P(x^{-1}) = x^{d-1} + \sum_{n=1}^{d-1} u_n
x^{d-n-1}.
\end{equation}

\begin{lemma}\label{lem1}
If $p(x)/q(x)$ is a convergent of $g_\vu(x)$ with the rate at least
$c$ then $P^*(x)p(x^d)/q(x^d)$ is also a convergent of $g_\vu(x)$
with the rate at least $cd-d+1$.
\end{lemma}

\proof We have
$$
\left|\left| g_\vu(x) - \frac{p(x)}{q(x)}\right|\right| \le
-2||q(x)|| -c.
$$
By substituting $x^d$ instead of $x$ and applying the functional
relation we get
$$
\left|\left| \frac{g_\vu(x)}{P^*(x)} -
\frac{p(x^d)}{q(x^d)}\right|\right| \le -2d||q(x)|| - cd.
$$
Multiply both sides of this equation by $P^*(x)$ and finally get
$$
\left|\left| g_\vu(x) - \frac{P^*(x)p(x^d)}{q(x^d)}\right|\right|
\le -2||q(x^d)|| - cd+d-1.
$$
\endproof

{\bf Remark.} This lemma shows the importance of the condition that
$||P^*(x)||<d$ and in turn of the condition $||P(x)||<d$. In this
case, any convergent $p(x)/q(x)$ with the rate of approximation
$c\ge 1$ allows us to construct an infinite series of other
convergents. Otherwise one needs the value of $c$ to be big enough,
so that $cd-||P^*(x)||>0$. However we can not guarantee that there
exists a convergent of $g_\vu(x)$ with the rate of approximation
strictly bigger than one.

By applying Lemma~\ref{lem1} several times we get the following
\begin{corollary}\label{cor1}
Let $k\in\NN$. If $p(x)/q(x)$ is a convergent of $g_\vu(x)$ with the
rate of approximation at least $c$ then
$$
\prod_{t=0}^{k-1} P^*(x^{d^t})\frac{p(x^{d^k})}{q(x^{d^k})}
$$
is also a convergent of $g_\vu(x)$ with the rate of approximation at
least $(c-1)d^k+1$.
\end{corollary}

Lemma~\ref{lem1} provides the following very nice criterium for
badly approximable series $g_\vu(x)$.

\begin{proposition}\label{lem2}
The series $g_\vu(x)$ is badly approximable if and only if all
partial quotients $a_n(x)$, $n\in\NN$ of its continued fraction are
linear.
\end{proposition}

\proof The ``if'' part of the lemma is straightforward. Let's show
the other part. Assume that $||a_n(x)||\ge 2$ for some $n\in\NN$.
Then the rate of approximation of convergent $r_1(x) :=
p_n(x)/q_n(x)$ is $c=c_1\ge 2$. Then by Lemma~\ref{lem1} there
exists another convergent $r_2(x)$ of $g(x)$ with the rate of
approximation $c_2\ge c_1d-d+1>c_1$. We use Lemma~\ref{lem1}
iteratively for convergents $r_2(x), r_3(x),...,r_m(x)$ to construct
a new convergent $r_{m+1}(x)$ with the rate of approximation
$c_{m+1}>c_m$. Hence $g_\vu(x)$ has a series of convergents with
unbounded rate of approximation which in turn implies that
$g_\vu(x)$ is well approximable. \endproof

With the help of Proposition~\ref{lem2} we show that the ``only if''
part of Theorem~\ref{th8} follows from Theorem~\ref{th5}. Indeed,
assume that $g_\vu(x)$ is not badly approximable. Then, by
Proposition~\ref{lem2}, there exists a partial quotient $a_{n_0}(x)$
of degree $c>1$. Then Theorem~\ref{th5} implies that, as soon as
$g_\vu(b)\neq 0$, $\mu(g_\vu(b))\ge 2+(c-1)/n_0 >2$.

In the rest of this section we look at the coordinates
$u_1,\ldots,u_{d-1}$ of $\vu$ as independent variables. Then the
Hankel determinant $H_n(f)$ of the series $f_\vu(x^{-1})$ is a
polynomial over them, i.e. $H_n(f) \in \FF[u_1,u_2,\ldots,
u_{d-1}]$.

\begin{lemma}\label{lem3}
Let $\chr\, \FF = 0$ and $d$ be a prime number. Then for any
$n\in\NN$, the polynomial $H_n(f)$ is not identically zero.
\end{lemma}
\proof To check this lemma we need to provide just one value of
$\vu$ (or respectively one polynomial $P(x)$) such that the series
$g_\vu(x)$ is badly approximable. That would imply by
Proposition~\ref{lem2} that all partial quotients of $g_\vu(x)$ are
linear and finally Theorem~H2 implies that values of $H_n(f)$ for
all $n\in\NN$ are non-zero, and therefore it is not zero
identically.

We use the technique which was firstly introduced by Han
in~\cite{han_2015a}. If $d=2$ then we know
from~\cite{badziahin_zorin_2014} that $g_{-1}(x)$ is badly
approximable. Let $d = p$ be an odd prime number. Then take
$$
P(x) = (x+1)^{\frac{p-1}{2}}.
$$
Consider the power series
$$
\tilde{f} (x) = f(x^{-1}) = \prod_{t=0}^\infty P(x^{d^t}).
$$
It satisfies the functional relation $\tilde{f}(x) = P(x)
\tilde{f}(x^p)$. Consider this equation over $\FF_p = \ZZ/p\ZZ$. It
becomes
$$
\tilde{f}(x)= (x+1)^{\frac{p-1}{2}} \tilde {f} (x)^p\quad\mbox{ or
}\quad \big(\tilde{f}(x)^2(x+1)\big)^{\frac{p-1}{2}} = 1,
$$
as $f(x)\not \equiv 0 \pmod p$. Therefore the series $\tilde{f}(x)$
is a solution of one of the equations $\tilde{f}(x)^2 (x+1) = a$
where $a$ is some quadratic residue over $\FF_p$. Definitely,
$\tilde{f}(x)$ can not be rational, therefore by~\cite[Theorem
1.1]{han_2015b} its Hankel continued fraction is ultimately periodic
which in turn yields that the sequence of non-zero values $H_n(f)$
over $\FF_p$ is also ultimately periodic. Going back to $\QQ$, the
Hankel continued fraction of $\tilde{f}$ is badly approximable and
hence $g(x)=x^{-1}\tilde{f}(x^{-1})$ is badly approximable.
\endproof

Lemma~\ref{lem3} only covers the case of prime $d$. Almost certainly
the same result should be true for any integer $d\ge 2$. It would be
interesting to see the proof of that statement. The author can
extend this lemma to integer powers of prime numbers, however the
other cases remain open.

We emphasize that the remaining results of this section are for
$\FF=\CC$ or for the subfields of $\CC$.

\begin{lemma}\label{lem4}
Let $d$ be prime. For any $\vu = (u_1,u_2,\ldots,
u_{d-1})\in\CC^{d-1}$ there exists a sequence of vectors $\vv_i\in
\CC^{d-1}$ such that $\vv_i \to \vu$ as $i\to\infty$ and all the
series $g_{\vv_i}(x)$ are badly approximable.
\end{lemma}

\proof Let $\vv\in\CC^{d-1}$. From Proposition~\ref{lem2} and
Theorem~H2 we know that $g_{\vv}(x)$ is badly approximable if and
only if all Hankel determinants $H_n(f_\vv)$ are not zero.
Lemma~\ref{lem3} implies that for every $n$ the equation
$H_n(f_\vv)=0$ is true for $\vv$ on a variety $\VVV_n$ of zero
Lebesgue measure. Whence, $g_\vv(x)$ is not badly approximable if
and only if $\vv$ belongs to countably many varieties with total
measure zero:
$$
\vv \in \RR^{d-1}\big/\bigcup_{n=1}^\infty \VVV_n.
$$

Take an arbitrary vector $\vu\in\RR^{d-1}$. For any $i\in\NN$ the
set
$$
S_i:= B(\vu, 1/i)\big/\bigcup_{n=1}^\infty \VVV_n
$$
is non-empty, where $B(\vu, r)$ is the ball in $\RR^{d-1}$ with the
center in $\vu$ and the radius $r$. Take any point $\vv_i\in S_i$.
By the construction $g_{\vv_i}(x)$ is badly approximable and also
$\vv_i\to \vu$ as $i\to \infty$. Hence the Lemma. \endproof

As we discussed before in Section~\ref{sec_contfrac}, for each
series $g_\vu(x)$ we associate partial quotients $a_n(x)$ and
parameters $\beta_n$ where $n\in \ZZ_{\ge 0}$. By
Proposition~\ref{lem2}, for badly approximable $g_\vu(x)$ all
polynomials $a_n(x)$ can be written as $a_n(x) = x+\alpha_n$.
Therefore we have a sequence of parameters $\beta_n$ and $\alpha_n$
which are uniquely defined by a badly approximable $g_\vu(x)$. It in
turn is defined by $\vu\in\CC^{d-1}$, hence we can look at
$\alpha_n$ and $\beta_n$ as maps:
$$
\alpha_n(\vu)\;:\; \CC^{d-1} \big/ \bigcup_{i=1}^n \VVV_i \to
\CC;\quad \beta_n(\vu)\;:\; \CC^{d-1} \big/ \bigcup_{i=1}^n
\VVV_n\to \CC;
$$

\begin{lemma}\label{lem5}
For each $n\in\ZZ_{\ge 0}$, the maps $\alpha_n(\vu)$ and
$\beta_n(\vu)$ are continuous.
\end{lemma}

\proof Firstly note that each coefficient $c_n$ in the formula
$$
g_\vu(x) = \sum_{n=1}^\infty c_nx^{-n}
$$
is a continuous function of $\vu$: $c_n(\vu)$.

Secondly, one can easily check that the $n$'th convergent
$p_n(x)/q_n(x)$ of badly approximable $g_\vu(x)$ is uniquely defined
by the first $2n+2$ terms of the series $g_\vu(x)$. Moreover, if
$q_n(x)$ is monic then $q_n(x)$ and $p_n(x)$ is a continuous map
from the coefficients $c_1,\ldots, c_{2n}$ to $\CC[x]$. Indeed, if
$q_n(x) = \sum_{i=0}^{n-1} a_ix^i+x^n$ then the coefficients can be
derived from the system
$$
\sum_{i=0}^{n-1} {c_{i+k}a_i} + c_{n+k} = 0
$$
for each $k$ between one and $n$. The matrix of this system is
basically $n$'th Hankel matrix which is invertible, because
$H_n(f_\vu)\neq 0$.

Finally, all terms $\alpha_i$ and $\beta_i$ are continuous maps from
$q_n(x)$ and $p_n(x)$ to $\CC$. The last statement follows from the
equation
$$
\frac{p_n(x)}{q_n(x)} = \bigk_{i=1}^n
\frac{\beta_{i+1}}{(x+\alpha_i)}.
$$ \endproof

\begin{lemma}\label{lem6}
Let $m\in \NN$ and $(\vu_i)_{i\in\NN}$ be the sequence of vectors in
$\CC^{d-1}$ with $\lim_{i\to \infty} \vu_i = \vu$ such that the
first $m$ partial quotients of $g_{\vu_i}(x)$ are linear. Assume
that for any $1\le n\le m$ there exist positive constants $c_n$ and
$C_n$ such that
\begin{equation}\label{eq_lem6}
\lim_{i\to \infty} |\alpha_n(\vu_i)|<C_n\;\mbox{ and
}\;c_n<\lim_{i\to \infty} |\beta_n(\vu_i)|<C_n.
\end{equation}
Then the first $m$ partial quotients of $g_{\vu}(x)$ are also linear
with coefficients
$$
\alpha_n(\vu) = \lim_{i\to \infty} \alpha_n(\vu_i),\quad
\beta_n(\vu) = \lim_{i\to \infty} \beta_n(\vu_i), \quad 1\le n\le m.
$$
\end{lemma}

The straightforward corollary of this lemma is that if
$g_{\vu_i}(x)$ are all badly approximable and \eqref{eq_lem6} is
satisfied for all $n\in\NN$ then the limiting series $g_\vu(x)$ is
also badly approximable.

\proof Let
$$
\frac{p_{n,\vu_i}(x)}{q_{n,\vu_i}(x)}
$$
be the $n$'th convergent of $g_{\vu_i}(x)$. Then we have
$$
\left|\left| g_{\vu_i}(x) -
\frac{p_{n,\vu_i}(x)}{q_{n,\vu_i}(x)}\right|\right| \le
-2||q_{n,\vu_i}(x)||-1
$$
Since $\alpha_n(\vu)$ and $\beta_n(\vu)$ are continuous, the limits
$\alpha_n = \lim_{i\to \infty} \alpha_n(\vu_i)$ and $\beta_n =
\lim_{i\to\infty}\beta_n(\vu_i)$ exist. From~\eqref{eq_lem6} we have
that $\alpha_n\le C_n$ and $0<c_n\le \beta_n\le C_n$. By continuity
we also have
$$
\frac{p_{n,\vu_i}(x)}{q_{n,\vu_i}(x)} = \bigk_{n=1}^m
\frac{\beta_n(\vu_i)}{x+\alpha_n(\vu_i)} \quad\to\quad \bigk_{n=0}^m
\frac{\beta_n}{x+\alpha_n} = \frac{p_n(x)}{q_n(x)}.
$$
Then again by continuity we have that the first $2||q_{n}(x)||+1$
terms of $g_{\vu_i}(x)$ tend to the corresponding terms of
$g_\vu(x)$. Therefore
$$
\left|\left| g_{\vu}(x) - \frac{p_{n}(x)}{q_{n}(x)}\right|\right|
\le -2||q_{n}(x)||-1,
$$
which in turn implies that $p_n(x)/q_n(x)$ are convergents of
$g_\vu(x)$.\endproof

\section{Irrationality exponents of
$g_\vu(b)$ for well approximable series}\label{sec_irr}

%
%
%
%
%

Throughout this section we assume that $g(x)$ is not badly
approximable. Proposition~\ref{lem2} asserts that in this case there
exists $n\in \NN$ such that the $n$-th convergent $p_n(x)/q_n(x)$
has rate of approximation $c\ge 2$. Then we can provide lower and
upper bounds for $\mu(b)$ which depend on the smallest value of $n$
with this property.

\smallskip\noindent\textsc{Proof of Theorem~\ref{th5}.} It is sufficient for any $\epsilon>0$ to provide an
infinite sequence of rational numbers $a_k/b_k$ such that
$$
\left|g_\vu(b) - \frac{a_k}{b_k}\right| <
\frac{\gamma}{b_k^{2+(c-1)n_0^{-1}-\epsilon}}.
$$
By construction of $n_0$ we have that $||q_n(x)||=n$ for all $n\le
n_0$ because all partial quotients of $g_\vu(x)$ are linear for
$n\le n_0$. Without loss of generality we may assume that both
$p_{n_0}(x)$ and $q_{n_0}(x)$ have integer coefficients. Indeed,
otherwise we just multiply both $p_{n_0}(x)$ and $q_{n_0}(x)$ by the
least common multiple of the denominators of all the coefficients of
both polynomials. We can also write $P^*(x)$ as $D^{-1}\tilde{P}(x)$
where $\tilde{P}(x)\in\ZZ[x]$ and $D\in\ZZ$.

Consider the following function
$$
F(x):= g_\vu(x) - \frac{p_{n_0}(x)}{q_{n_0}(x)}.
$$
It can be written as an infinite series and moreover, since
$p_{n_0}(x)/q_{n_0}(x)$ is a convergent of $g_\vu(x)$ with rate of
approximation $c$, we have
$$
F(x) = \sum_{n=2n_0+c}^\infty c_n x^{-n}\quad c_n\in\RR.
$$
We know that $F(x)$ converges absolutely for all $|x|>1$ and
therefore for all $|x|\ge 2$ we have
\begin{equation}\label{eq_f}
|x^{2n_0+c}F(x)| =\left|\sum_{n=0}^\infty c_{n+2n_0+c}x^{-n}\right|
\le \sum_{n=0}^\infty |c_{n+2n_0+c}|2^{-n} =:\gamma_1.
\end{equation}
In other words there exists an absolute constant $\gamma_1$ such
that for all $|x|\ge 2$ we have $|F(x)|\le \gamma_1 x^{-2n_0-c}$.

Now apply the functional equation~\eqref{eq_g} $k$ times to get
\begin{equation}\label{eq_th5}
F(x^{d^k})\prod_{t=0}^{k-1} P^*(x^{d^t}) = g_\vu(x) -
\frac{p_{n_0}(x^{d^k})\prod_{t=0}^{k-1}P^*(x^{d^t})}{q_{n_0}(x^{d^k})}
= g_\vu(x) -
\frac{p_{n_0}(x^{d^k})\prod_{t=0}^{k-1}\tilde{P}(x^{d^t})}{D^kq_{n_0}(x^{d^k})}.
\end{equation}
We set $a_k:= p_{n_0}(b^{d^k})\prod_{t=0}^{k-1}\tilde{P}(b^{d^t})$
and $b_k:=D^kq_{n_0}(b^{d^k})$. By construction, they are both
integer. Moreover, one can check that
$$
\lim_{k\to \infty} \frac{q_{n_0}(b^{d^k})}{b^{n_0d^k}} = \gamma_2,
$$
where $\gamma_2$ is the leading coefficient of $q_{n_0}(x)$.
Therefore for large enough $k$ we have $|b_k|\le 2|\gamma_2|D^k
b^{n_0d^k}$.

Now we use inequality~\eqref{eq_f} for $F(x)$ and~\eqref{eq_th5} to
estimate $|g_\vu(b) - a_k/b_k|$:
$$
\left|g_\vu(b) - \frac{a_k}{b_k}\right| \le \frac{\gamma_1
\prod_{t=0}^{k-1} P^*(b^{d^t})}{b^{(2n_0+c)d^k}} = \frac{\gamma_1
\prod_{t=0}^{k-1} P(b^{-d^t})}{b^{(2n_0+c-1)d^k+1}}.
$$
Since $\prod_{t=0}^{\infty} P(x^{-d^t})$ converges absolutely for
all $|x|>1$, there exists a uniform upper bound $\gamma_3$ such that
for all $|b|\ge 2$ we have $|\prod_{t=0}^{\infty} P(b^{-d^t})|\le
\gamma_3$. Next, by solving the equation
$$
(D^kb^{n_0d^k})^x = b^{(2n_0+c-1)d^k}
$$
we get
$$
x =x(k) =  \left(2+\frac{c-1}{n_0}\right)\left(1+\frac{k\log_b
D}{n_0d^k}\right)^{-1}
$$
As $k$ tends to infinity, $x(k)$ tends to $2+(c-1)/n_0$. Therefore
for any $\epsilon>0$ we can find $k(\epsilon)$ large enough so that
for any $k>k(\epsilon)$,  $x(k) > 2+(c-1)/n_0 - \epsilon$ and
therefore
$$
\left|g_\vu(b) - \frac{a_k}{b_k}\right| <
\frac{(2\gamma_2)^{2+(c-1)/n_0}\gamma_1\gamma_3}{b} \cdot
b_k^{-2-(c-1)/n_0 + \epsilon}.
$$
\endproof

\section{Recurrent formulae for continued fractions of
$g_\vu(x)$}\label{sec_conds}

In this section we construct the continued fraction of the series
$g_\vu(x)$. Throughout the whole section we assume that $g_\vu(x)$
is badly approximable. Then, by Proposition~\ref{lem2}, its
continued fraction is determined by the terms $\alpha_n$ and
$\beta_n$ where the partial quotients $a_n(x) = x+\alpha_n$ and the
parameters $\beta_n$ satisfy the recurrent
formulae~\eqref{mconv_def}.

\begin{proposition}\label{prop1}
Let $g_\vu(x)$ be badly approximable. Then for any $k\in\ZZ_{\ge 0}$
one has
\begin{equation}\label{prop1_eq}
K_d^1 (\ova_{dk+1,d(k+1)},\ovb_{dk+1,d(k+1)}) = \beta_{dk+1}P^*(x),
\end{equation}
where $P^*(x)$ is given in~\eqref{pstar_def}.
\end{proposition}

\proof

Let $p_k(x)/q_k(x)$ be $k$th convergent of $g_\vu(x)$.
Proposition~\ref{lem2} asserts that $||q_k(x)|| = k$. We know from
Lemma~\ref{lem1} that $P^*(x)p_k(x^d)/q_k(x^d)$ is another
convergent of $g(x)$. We can assume that $P^*(x)p_k(x)$ and $q_k(x)$
are coprime. Indeed otherwise one can cancel their common divisor
from the fraction $P^*(x)p_k(x^d)/q_k(x^d)$ and its rate of
convergence will become bigger than one, which contradicts to
Proposition~\ref{lem2}. Therefore we get that the fraction
$P^*(x)p_k(x^d)/q_k(x^d)$ is in fact $dk$'th convergent of
$g_\vu(x)$. By following this arguments for each $k\in\NN$ we get
that
$$
p_{dk}(x) \equiv 0\pmod{P^*(x)}.
$$

Consider the equation~\eqref{lem7_eq} from Lemma~\ref{lem7} with
$n=dk$ and $m=d$ modulo $P^*(x)$:
$$
0\equiv p_{d(k+1)}(x) \equiv
K_d^1(\ova_{dk+1,d(k+1)},\ovb_{dk+1,d(k+1)})p_{dk-1}(x)\pmod{P^*(x)}.
$$
The observation $\gcd(p_{dk-1}(x),p_{dk}(x)) =
\gcd(p_{dk-1}(x),P^*(x)) = 1$ implies that
$$
K_d^1(\ova_{dk+1,d(k+1)},\ovb_{dk+1,d(k+1)})\equiv 0\pmod{P^*(x)}.
$$
By~\eqref{cont_deg} the degree of the left hand side coincides with
those of $P^*(x)$. Then comparing the leading coefficients of the
polynomials in the congruence finishes the proof.
\endproof

Polynomial equation~\eqref{prop1_eq} gives us $d-1$ relations
between various values $\alpha_n$ and $\beta_n$ for each
$k\in\ZZ_{\ge 0}$. We just need to compare the corresponding
coefficients of the polynomials from both sides of the equation.
However they are still not enough to provide the recurrent formula
for all values $\alpha_{dk+1},\ldots, \alpha_{d(k+1)},
\beta_{dk+1},\ldots, \beta_{d(k+1)}$. More relations can be derived
from the following:

\begin{proposition}\label{prop2}
Let $g_\vu(x)$ be badly approximable. Then for any $k\in\NN$ one has
\begin{equation}\label{prop2_eq1}
K_{2d}^1(\ova_{dk+1,d(k+2)},\ovb_{dk+1,d(k+2)}) =
\beta_{dk+1}(x^d+\alpha_{k+2})P^*(x)
\end{equation}
and
\begin{equation}\label{prop2_eq2}
K_{2d}^0(\ova_{dk+1,d(k+2)},\ovb_{dk+1,d(k+2)}) =
\beta_{k+2}+(x^d+\alpha_{k+2})K_d^0(\ova_{kd+1,k(d+1)},\ovb_{kd+1,k(d+1)}).
\end{equation}
\end{proposition}

\proof

For convenience we will use the following notation throughout the
proof: $K_{2d}^0(x) :=
K_{2d}^0(\ova_{dk+1,d(k+2)},\ovb_{dk+1,d(k+2)})$,
$K_d^0(x):=K_d^0(\ova_{kd+1,d(k+1)},\ovb_{kd+1,d(k+1)})$. The
notions of $K_{2d}^1(x)$ and $K_d^1(x)$ are defined by analogy.

We provide two different relations between $q_{kd}(x),
q_{(k+1)d}(x)$ and $q_{(k+2)d}(x)$. The first one comes from the
fact that for each $m\in \NN$, $q_{md}(x) = q_{m}(x^d)$, which was
shown in the proof of Proposition~\ref{prop1}. Therefore the
application of~\eqref{mconv_def} gives us
\begin{equation}\label{eq1_prop2}
q_{(k+2)d}(x) = (x^d + \alpha_{k+2})q_{(k+1)d}(x) +
\beta_{k+2}q_{kd}(x).
\end{equation}

On the other hand~\eqref{lem7_eq} implies
$$
\begin{array}{rl}
q_{(k+1)d}(x) =& K_d^0(x) q_{kd}(x)+ K_d^1(x)q_{kd-1}(x)\\[1ex]
\ [\mbox{by Proposition~\ref{prop1}}] =& K_d^0(x) q_{kd}(x)+
\beta_{kd+1}P^*(x)q_{kd-1}(x).
\end{array}
$$
From this formula we can write $q_{kd-1}(x)$ in terms of $q_{kd}(x)$
and $q_{(k+1)d}(x)$.
$$
q_{kd-1}(x) = \frac{q_{(k+1)d}(x) -
K_d^0(x)q_{kd}(x)}{\beta_{kd+1}P^*(x)}.
$$
Next,~\eqref{lem7_eq} also gives us
$$
\begin{array}{rl}
q_{(k+2)d}(x) =& K_{2d}^0(x)q_{kd}(x) +
K_{2d}^1(x)q_{kd-1}(x)\\[1ex]
=&\displaystyle\left(K_{2d}^0(x) -
\frac{K_{2d}^1(x)K_d^0(x)}{\beta_{kd+1}P^*(x)}\right)q_{kd}(x) +
\frac{K_{2d}^1(x)}{\beta_{kd+1}P^*(x)}q_{(k+1)d}(x).
\end{array}
$$

Combining the last formula with~\eqref{eq1_prop2} gives
\begin{equation}\label{eq2_prop2}
\left(\beta_{k+2} + \frac{K_{2d}^1(x)K_d^0(x)}{\beta_{kd+1}P^*(x)} -
K_{2d}^0(x)\right) q_{kd}(x) =
\left(\frac{K_{2d}^1(x)}{\beta_{kd+1}P^*(x)} - x^d-\alpha_{k+2}
\right)q_{(k+1)d}(x).
\end{equation}

Adapting the formula~\eqref{lem7_eq} to $K_{2d}^1$ gives
\begin{equation}\label{eq3_prop2}
K_{2d}^1(x) \!=\!
K_d^0(\ova_{kd+d+1,(k+2)d},\ovb_{kd+d+1,(k+2)d})K_d^1(x) +
K_d^1(\ova_{kd+d+1,(k+2)d},\ovb_{kd+d+1,(k+2)d})K_{d-1}^1(x)
\end{equation}
By Proposition~\ref{prop1} we get that $K_d^1(x)=\beta_{kd+1}P^*(x)$
and $K_d^1(\ova_{kd+d+1,(k+2)d},\ovb_{kd+d+1,(k+2)d}) =
\beta_{kd+d+1}P^*(x)$. This straightforwardly implies that the
expressions on the left and right hand sides of~\eqref{eq2_prop2}
are in fact polynomials. Moreover, since the leading coefficient of
$K_{2d}^1(x)$ is $\beta_{kd+1}$ we have that the degree of the
polynomial
$$
D(x):= \frac{K_{2d}^1(x)}{\beta_{kd+1}P^*(x)} - x^d-\alpha_{k+2}
$$
is at most $d-1$.

Two polynomials $q_{kd}(x)=q_k(x^d)$ and
$q_{(k+1)d}(x)=q_{k+1}(x^d)$ are coprime. Therefore $D(x)$ should be
a multiple of $q_{kd}(x)$. However for $k\ge 1$ its degree
$||q_{kd}(x)||$ is strictly bigger than $||D(x)||$ which is only
possible when $D(x) = 0$. This immediately gives the
formula~\eqref{prop2_eq1}. Finally,~\eqref{prop2_eq2} can be
achieved by equating the right hand side of~\eqref{eq2_prop2} to
zero. \endproof

\subsection{Recurrent formulae for small $d$}

Relations from Proposition~\ref{prop1} and~\ref{prop2} appear to be
enough to provide the recurrent formulae for the values $\alpha_n$
and $\beta_n$. We demonstrate that by constructing the recurrent
formulae for small values of $d$.

\medskip
\noindent{\bf The case $d=2$.} We have
$$
\vu = u,\quad g_u(x) = \prod_{t=0}^\infty (1+ux^{-d^t})\;\;\mbox{
and }\;\;P^*(x) = x+u.
$$

\smallskip\noindent\textsc{Proof of Theorem~\ref{th1}}. By Proposition~\ref{prop1} we have that for any $k\ge 0$,
$$
K_2^1(\ova_{2k+1,2k+2},\ovb_{2k+1,2k+2}) =
\beta_{2k+1}(x+\alpha_{2k+2}) = \beta_{2k+1} (x+u).
$$
Since we assumed that $g_u(x)$ is badly approximable,
$\beta_{2k+1}\neq 0$ and the formula straightforwardly implies that
$\alpha_{2k+2} = u$ for any $k\in\ZZ_{\ge 0}$.

Then we apply Proposition~\ref{prop2}. From~\eqref{prop2_eq1} for
any $k\in\NN$ we have
$$
\begin{array}{rl}
&((x+\alpha_{2k+4})(x+\alpha_{2k+3})(x+\alpha_{2k+2}) +
(x+\alpha_{2k+4})\beta_{2k+3} +
(x+\alpha_{2k+2})\beta_{2k+4})\beta_{2k+1}\\
=& \beta_{2k+1}(x^2+\alpha_{k+2})(x+u).
\end{array}
$$
We already know that $\alpha_{2k+2} = \alpha_{2k+4} = u$. Then
comparing the coefficients for $x^2, x$ and 1 gives
$$
\alpha_{2k+3} = -u;\quad \beta_{2k+3}+\beta_{2k+4} =
\alpha_{k+2}+u^2.
$$
Finally, look at equation~\eqref{prop2_eq2} modulo $K_2^0(x):=
K_2^0(\ova_{2k+1,2k+2},\ovb_{2k+1,2k+2})$:
\begin{equation}\label{eq_d2}
K_4^0(\ova_{2k+1,2k+4},\ovb_{2k+1,2k+4}) \equiv
(x\!+\!\alpha_{2k+1})(x\!+\!\alpha_{2k+4})\beta_{2k+3} =
(x\!+\!\alpha_{2k+1})(x\!+\!\alpha_{2k+2})\beta_{2k+3}.
\end{equation}
The right hand side of~\eqref{prop2_eq2} is congruent to
$\beta_{k+2}$ modulo $K_2^0(x)$. We also have,
$$
K_2^0(x) = (x+\alpha_{2k+2})(x+\alpha_{2k+1}) + \beta_{2k+2}
$$
and therefore the last expression in~\eqref{eq_d2} is congruent to
$-\beta_{2k+2}\beta_{2k+3}$. Hence this provides the following
relation between $\beta$'s:
$$
\beta_{2k+2}\beta_{2k+3} = -\beta_{k+2}.
$$
We collect all the data together and get the recurrent formulae
which allow us to confirm formulae~\eqref{recur_d2} for $\alpha_n$
and $\beta_n$ starting from $n=5$: for any $k\ge 1$,
$$
\begin{array}{l}
\alpha_{2k+2} = u;\quad \alpha_{2k+3} = -u;\\[1ex]
\beta_{2k+3} = -\frac{\beta_{k+2}}{\beta_{2k+2}},\quad \beta_{2k+4}
= \alpha_{k+2}+u^2 - \beta_{2k+3}.
\end{array}
$$

To finish the proof we need to find the values
$\alpha_1,\ldots,\alpha_4$ and $\beta_1,\ldots, \beta_4$. By direct
computation one can easily check that the first convergent of
$g_u(x)$ is $(x-u)^{-1}$. That together with Lemma~\ref{lem1} gives
us
$$
\frac{p_1(x)}{q_1(x)} = \frac{1}{x-u};\quad \frac{p_2(x)}{q_2(x)} =
\frac{x+u}{x^2 - u},\quad \frac{p_4(x)}{q_4(x)} =
\frac{(x+u)(x^2+u)}{x^4-u}.
$$
We find the denominator $q_3(x) = x^3+ax^2+bx+c$ of the third
convergent by noticing that
$$
g_u(x) = x^{-1} + ux^{-2}+ux^{-3}+u^2x^{-4} + ux^{-5} + u^2x^{-6}
+\ldots
$$
and that the coefficients for $x^{-1}, x^{-2}$ and $x^{-3}$ of the
expression $ g_u(x)(x^3 + ax^2+bx+c)$ are all zeroes. That gives us
the system of linear equations in $a,b,c$ with solutions $a=-u,
b=-u-1, c=u(u+1)$. That finally gives us
$$
\frac{p_3(x)}{q_3(x)} = \frac{x^2 - u^2-1}{(x-u)(x^2 - u - 1)}
$$
These convergents give us the initial values:
$$
\alpha_1 = -u,\; \alpha_2 = u,\; \alpha_3 = -u,\; \alpha_4 = u;
$$$$
\beta_1 = 1,\; \beta_2 = u^2 - u,\; \beta_3 = -1,\; \beta_4 =
u^2+u+1 .
$$
Now we have all the relations from~\eqref{recur_d2}. \endproof

%

\medskip\noindent{\bf The case d = 3.} We have
$$
\vu = (u,v),\quad g_\vu(x) = \prod_{t=0}^\infty
(1+ux^{-3^t}+vx^{-2\cdot 3^t})\quad\mbox{and}\quad P^*(x) =
x^2+ux+v.
$$

\smallskip\noindent\textsc{Proof of Theorem~\ref{th2}}. We proceed in a similar way as for the case $d=2$.
Proposition~\ref{prop1} gives us that for any integer $k\ge 0$,
\begin{equation}\label{eq1_d3}
(x+\alpha_{3k+2})(x+\alpha_{3k+3}) + \beta_{3k+3} = x^2+ux+v
\end{equation}
This immediately implies some relations between the coefficients:
\begin{equation}\label{rel_d3}
\alpha_{3k+2}+\alpha_{3k+3} = u,\quad
\alpha_{3k+2}\alpha_{3k+3}+\beta_{3k+3} = v.
\end{equation}

Next, we apply the equation~\eqref{prop2_eq1}, where we
use~\eqref{eq3_prop2} to compute $K_{2d}^1(x)$. For $k\ge 1$ we get
$$
\begin{array}{rl}
&\beta_{3k+1}(x^2\!+\!ux\!+\!v)((x+\alpha_{3k+4})(x^2\!+\!ux\!+\!v)+(x\!+\!\alpha_{3k+6})\beta_{3k+5}+(x\!+\!\alpha_{3k+2})\beta_{3k+4})\\[1ex]
=\!\!\!\!&\beta_{3k+1}(x^2+ux+v) (x^3+\alpha_{k+2}).
\end{array}
$$
Comparing the coefficients then gives
\begin{equation}\label{rel2_d3}
\begin{array}{l}
\alpha_{3k+4} + u = 0,\quad
u\alpha_{3k+4}+v+\beta_{3k+5}+\beta_{3k+4} = 0,\\[1ex]
\alpha_{3k+4}v + \alpha_{3k+6}\beta_{3k+5} +
\alpha_{3k+2}\beta_{3k+4} = \alpha_{k+2}.
\end{array}
\end{equation}

Finally, as before, apply the equation~\eqref{prop2_eq2} modulo
$K_3^0(x)$:
$$
K_3^0(x)= K_3^0(\ova_{3k+1,3k+3},\ovb_{3k+1,3k+3}) =
(x+\alpha_{3k+1})(x^2+ux+v)+(x+\alpha_{3k+3})\beta_{3k+2}.
$$
We get
$$
\begin{array}{rl}
\beta_{k+2}\equiv&
\beta_4((x+\alpha_{3k+5})(x+\alpha_{3k+6})+\beta_{3k+6})((x+\alpha_{3k+2})(x+\alpha_{3k+1})+\beta_{3k+2})\\[1ex]
\stackrel{\eqref{eq1_d3}}{\equiv} &\beta_{3k+4}
P^*(x)((x+\alpha_{3k+2})(x+\alpha_{3k+1})+\beta_{3k+2})\\[1ex]
\equiv& \beta_{3k+4}( (x+\alpha_{3k+2})K_3^0(x) -
(x+\alpha_{3k+3})(x+\alpha_{3k+2})\beta_{3k+2} +
P^*(x)\beta_{3k+2})\\[1ex]
\stackrel{\eqref{eq1_d3}}{\equiv}&
\beta_{3k+2}\beta_{3k+3}\beta_{3k+4}
\end{array}
$$
or $\beta_{k+2} = \beta_{3k+2}\beta_{3k+3}\beta_{3k+4}$. Combining
this formula with~\eqref{rel_d3} and~\eqref{rel2_d3} we finally get
recurrent formulae for all values $\alpha_n$ and $\beta_n$
satisfy~\eqref{recur_d3} starting from $n=7$:
$$
\begin{array}{l}
\displaystyle \alpha_{3k+4} = -u,\quad \beta_{3k+4} =
\frac{\beta_{k+2}}{\beta_{3k+3}\beta_{3k+2}};\\[1ex]
\displaystyle \beta_{3k+5} = u^2 - v - \beta_{3k+4},\quad
\alpha_{3k+5} = u -
\frac{\alpha_{k+2}+uv-\alpha_{3k+2}\beta_{3k+4}}{\beta_{3k+5}}\\[1ex]
\alpha_{3k+6} = u-\alpha_{3k+5},\quad \beta_{3k+6} = v -
\alpha_{3k+5}\alpha_{3k+6}.
\end{array}
$$
Note that since~\eqref{rel_d3} is also true for $k=0$, $\alpha_6$
and $\beta_6$ can also be computed by~\eqref{recur_d3}. Therefore it
remains to compute $\alpha_1,\ldots,\alpha_5$ and $\beta_1,\ldots
\beta_5$. We do that straight from calculating the first five
convergents of $g_\vu(x)$. To save the space we will only provide
their denominators.
$$
q_1(x) = (x-u),\quad q_2(x) = x^2 + \frac{u(v-1)}{v-u^2}x +
\frac{u^2 - v^2}{v-u^2},\quad q_3(x) = x^3 - u
$$
This allows us to get the values:
$$
\alpha_1 = -u, \alpha_2 = \frac{u(2v-1-u^2)}{v-u^2},\alpha_3 =
\frac{-u(v-1)}{v-u^2};
$$
$$
\beta_1 = 1, \beta_2 = u^2-v,\beta_3 = \frac{u^2+u^4+v^3 -
3u^2v}{(v-u^2)^2}.
$$
Finally we use Mathematica to compute
$\alpha_4,\alpha_5,\beta_4,\beta_5$ and to confirm that they satisfy
the recurrent equations~\eqref{recur_d3} with $k=0$. This gives us
all relations from~\eqref{init_d3} and~\eqref{recur_d3}. \endproof

\section{Badly approximable series $g_\vu(x)$ for small
$d$.}\label{subsec_bad23}

Theorems~\ref{th1} and~\ref{th2} are only valid for badly
approximable series $g_\vu(x)$. In this section we try to answer the
question: for what values $\vu$ the series $g_\vu(x)$ is in fact
badly approximable? Then for such series $g_\vu(x)$ all machinery of
the previous paragraph can be used.

\smallskip\noindent\textsc{Proof of Theorem \ref{th7}.} Assume that the first $n$
terms $\alpha_k$ and $\beta_k$ satisfy
\begin{equation}\label{eq_th7}
|\alpha_k|<\infty,\; 0<|\beta_k|<\infty;\quad 1\le k\le n.
\end{equation}
Then by Lemma~\ref{lem4} there exists a sequence of vectors $\vu_i$
such that $\vu_i\to \vu$ and $g_{\vu_i}(x)$ are badly approximable.
Their parameters $\alpha_k(\vu_i)$ and $\beta_k(\vu_i)$ are computed
by formulae~\eqref{recur_d2} for $d=2$ and by~\eqref{init_d3}
and~\eqref{recur_d3} for $d=3$. Therefore by continuity the
coefficients $\alpha_k(\vu_i)$ and $\beta_k(\vu_i)$ tend to
$\alpha_k$ and $\beta_k$ respectively. Moreover, for each $k$ there
exists $i_0= i_0(k)$ and $c_k = |b_k|/2, C_k = \max\{2|\alpha_k|,
2|\beta_k|\}$ such that for $i>i_0$ we have $\alpha_k(\vu_i)<C_k$
and $c_k<\beta_k(\vu_i)<C_k$. Finally, Lemma~\ref{lem7} confirms
that $\alpha_k$ and $\beta_k$ with $1\le k\le n$ are indeed the
coefficients of the continued fraction of $g_\vu(x)$.

Note that all division in formulae~\eqref{recur_d2},~\eqref{init_d3}
and~\eqref{recur_d3} for $\alpha_k$ and $\beta_k$ are by some values
of $\beta_m$ with $m\le k$. Therefore, as soon as $\beta_1,\ldots,
\beta_n$ do not vanish, the condition~\eqref{eq_th7} is
automatically satisfied. Moreover, in this case we also can not have
$|\beta_{n+1}|=\infty$.

Finally, assume that $\beta_{n+1}=0$. If the $(n+1)$th partial
quotient of $g_\vu(x)$ is linear than by Lemma~\ref{lem5} the
sequence $\beta_{n+1}(\vu_i)$ should tend to $\beta_{n+1}(\vu)$.
Therefore $\beta_{n+1}(\vu)=0$ which is impossible, because all
values of $\beta_k$ in a continued fraction for $g_\vu(x)$ must be
non-zero. Hence we have a contradiction and the $(n+1)$th partial
quotient in $g_{\vu}(x)$ is not zero. \endproof

\noindent{\bf The case $d=2$.}

\smallskip\noindent\textsc{Proof of Theorem~\ref{th3}.} By
Theorem~\ref{th7}, $g_\vu(x)$ is well approximable if and only if
one of the values $\beta_n$ vanishes. From~\eqref{recur_d2} there
are two obvious values $u=0$ and $u=1$ when $\beta_2 = 0$. They in
fact produce rational functions: $g_0(x) = x^{-1}, g_1(x) =
(x-1)^{-1}$. On the other hand the values of $\beta_3$ and $\beta_4$
do not equal to zero for any rational (and in fact any real) values
of $u$.

\begin{lemma}\label{lem8}
For any $n\ge 3$ the value $\beta_n(u)$ can be written as
$$
\beta_n(u) = \frac{e_n(u)}{d_n(u)},
$$
where $e_n(u), d_n(u)\in\ZZ[x]$ and the leading and constant
coefficients of both $e_n(u), d_n(u)$ equal $\pm 1$.
\end{lemma}
\proof It can be easily checked by induction. It is true for $n=3$
and $n=4$. We assume that the statement is true for all values
$\beta_3(u), \ldots \beta_{2k+2}(u)$ and prove it for
$\beta_{2k+3}(u)$ and $\beta_{2k+4}(u)$. In addition we will check
the following condition: $||e_{2k}(u)|| = ||d_{2k}(u)||+2$,
$||e_{2k+1}(u)||\le ||d_{2k+1}(u)||$.

By~\eqref{recur_d2} we have
$$
\beta_{2k+3}(u) =
\frac{-d_{2k+2}(u)e_{k+2}(u)}{e_{2k+2}(u)d_{k+2}(u)}.
$$
Its numerator and denominator clearly satisfy the conditions of the
lemma together with
$$
\begin{array}{rl}
||e_{2k+3}(u)|| &=  ||d_{2k+2}(u)e_{k+2}(u)|| = ||e_{2k+2}(u)|| +
||e_{k+2}(u)|| - 2\\
&\le ||e_{2k+2}(u)d_{k+2}(u)|| = ||d_{2k+3}(u)||.
\end{array}
$$
For $\beta_{2k+4}(u)$ we have
$$
\beta_{2k+4}(u) = \alpha_{k+2}(u) + u^2 -
\frac{e_{2k+3}(u)}{d_{2k+3}(u)} = \frac{(u^2 \pm u)d_{2k+3}(u) -
e_{2k+3}(u)}{d_{2k+3}(u)}.
$$
The leading coefficient of the numerator on the right hand side
comes from $u^2d_{2k+3}(u)$ and the constant coefficient comes from
$-e_{2k+3}(u)$. Both of them by inductional hypothesis are plus or
minus one. Since $||e_{2k+3}(u)||\le ||d_{2k+3}(u)||$, we have
$$
||e_{2k+4}(u)||=||(u^2 \pm u)d_{2k+3}(u) - e_{2k+3}(u)|| =
||u^2d_{2k+3}(u)|| = 2+||d_{2k+3}(u)||=||d_{2k+4}(u)||.
$$
This completes the induction.
\endproof

The obvious corollary from Lemma~\ref{lem8} is that if $u\in\QQ$ and
$g_n(u) = 0$ for some $n\ge 3$ then $u$ is either plus or minus one.
Indeed these are the only possible rational roots of the equation
$e_n(u) = 0$. On the other hand it was shown
in~\cite{badziahin_zorin_2015} that $g_{-1}(x)$ is badly
approximable. Theorem~\ref{th3} is proven.

\noindent{\bf Remark.} There exist real values of $u$ such that
$g_u(x)$ is well approximable. For example, one can check that if
$u$ is any real root of the equation $u^4 - u - 1=0$ then
$\beta_6(u) = 0$.

\medskip\noindent{\bf The case $d=3$.}

\noindent\textsc{Proof of Theorem~\ref{th9}}. As for $d=2$ we
investigate the case when $\beta_n=0$ for some $n$.
From~\eqref{init_d3}, the equation $\beta_2 = 0$ gives an infinite
collection of vectors $\vu = (u,u^2)$ such that the series:
$$
g_{(u,u^2)} (x) = x^{-1}\prod_{t=0}^\infty (1 + ux^{-3^t} + u^2
x^{-2\cdot 3^t})
$$
is well approximable. Theorem~\ref{th5} with $n_0=1$ and $c\ge 2$ then
asserts that for any $u\in\QQ$ and integer $b\ge 2$, as soon as $g_{(
u,u^2)}(x)\not\in\QQ(x)$ and $g_{(u,u^2)}(b)\neq 0$, we have
$\mu(g_{(u,u^2)}(b))\ge 3$.

The equation $\beta_3=0$ can be written as $u^2+u^4+v^3-3u^2v=0$. It
gives an infinite parametrised series of rational solutions: $u=
s^3$, $v = -s^2(s^2+1)$ where $s\in\QQ$. It has only one
intersection with the collection $(u, u^2)$ above, namely when
$s=0$. This solution can be ignored, because $g_{(0,0)}(x)=x^{-1}$
is a rational function. Hence we have another set of well
approximable series:
$$
q_{(s^3,-s^2(s^2+1))}(x) = x^{-1}\prod_{t=0}^\infty (1+
s^3x^{-3^t}-s^2(s^2+1)x^{-2\cdot 3^t}),\quad s\in\QQ\backslash\{0\}.
$$
Direct computation shows that the second convergent of
$g_{(s^3,-s^2(s^2+1))}(x)$ is
$$
\frac{p_2(x)}{q_2(x)}= \frac{x+s(s^2+1)}{x^2+sx+s^2}
$$
and
$$
(x^2+sx+s^2)x^{-1}\prod_{t=0}^\infty
(1+s^3x^{-3^t}-s^2(s^2+1)x^{-2\cdot 3^t}) = -(s^6+s^4+s^2)x^{-5} +
\ldots
$$
If $s\neq 0$ the term $-(s^6+s^4+s^2)$ is non-zero and therefore
$$
\left|\left|g_{s^3,-s^2(s^2+1)}(x) -
\frac{p_2(x)}{q_2(x)}\right|\right| = -7
$$
or in other words the rate of approximation of the second convergent is
three. Then the application of Theorem~\ref{th5} with $n_0=2$ and $c=3$ tells
us that if $g_{s^3,-s^2(s^2+1)}(x)\not\in\QQ(x)$ and
$g_{s^3,-s^2(s^2+1)}(b)\neq 0$, we have $\mu(g_{s^3,-s^2(s^2+1)}(b))\ge
2+\frac22=3$ for all $s\in\QQ$ and all integer $b\ge 2$.

There is at least one less trivial example of well approximable
series. One can note that $g_{(2,1)}(x)$ is well approximable by
noticing that $\beta_6(2,1) = 0$. Direct computation shows that for
the fifth convergent of $g_{(2,1)}(x)$ we have $p_5(x) =
x^4+x^3+2x^2+4$, $q_5(x) = x^5-x^4+x^3-x^2+x-1$ and
$$
(x^5-x^4+x^3-x^2+x-1)x^{-1}\prod_{t=0}^\infty
(1+2x^{-3^t}+x^{-2\cdot 3^t}) = 3x^{-8}+\ldots.
$$
Therefore the rate of approximation of the fifth convergent of $g_{(2,1)}(x)$
is three. Since for any $|b|>1$ the value $g_{(2,1)}(b)$ is non-zero,
Theorem~\ref{th5} with $n_0=5$ and $c=3$ implies that $\mu(g_{(2,1)}(b))\ge
12/5$.\endproof

\smallskip\noindent\textsc{Proof of theorem~\ref{th4}}.
One can notice that $ g_{(0,v)}(x) = xg_{(v,0)}(x^2)$ and therefore
$g_{(0,v)}(x)$ is badly approximable if and only if so is
$g_{(v,0)}(x)$. Therefore without loss of generality we can only
assume the case $u=0$.

Let $u=0$. Then formulae~\eqref{init_d3},~\eqref{recur_d3} and an
easy induction give us that $\alpha_k=0$ and $\beta_{3k+3}=v$ for
all $k\ge 0$. We can write $\beta_k(\vu)$ as a rational function of
$v$:
$$
\beta_k(v) =: \frac{e_k(v)}{d_k(v)}
$$
where $e_k(v)$ and $d_k(v)$ are polynomials with integer
coefficients.

\begin{lemma}\label{lem9}
If $u=0$ then for any $k\in \NN$ the values $\beta_k$ satisfy the
following properties:
\begin{enumerate}
\item $||\beta_{3k+1}(v)||\le -1$ and
$||\beta_{3k+2}(v)||=||\beta_{3k+3}(v)||=1$;
\item The leading coefficient of $e_k(v)$ as well as of $d_k(v)$ is either plus or minus one;
\item If $|v|\ge 3$ then $|\beta_{3k+1}(v)|<1$ and
$|v|-1<|\beta_{3k+2}(v)|<|v|+1$.
\end{enumerate}
\end{lemma}

\proof All these items can simultaneously be shown by induction. For
$k=1$ one can easily check that:
$$
\beta_4(v) = \frac1v;\quad \beta_5(v) = -\frac{v^2+1}{v};\quad
\beta_6(v)=v.
$$
Also $\beta_{3k+3} = v$ obviously satisfies all the conditions for
each $k\in\NN$.

Assume that the properties are true for all integer values up to $k$
and prove it for $k+1$. By~\eqref{recur_d3} we have that
$$
||\beta_{3k+4}(v)|| =
||\beta_{k+2}(v)||-||\beta_{3k+2}(v)||-||\beta_{3k+3}(v)||\le -1;
$$
$$
|\beta_{3k+4}(v)| \le \frac{|v|+1}{|v|(|v|-1)}.
$$
The last expression is always less than one for $|v|\ge 3$. Next,
since we already know that $||\beta_{3k+4}(v)||<0$,
$$
||\beta_{3k+5}(v)|| = ||-v-\beta_{3k+4}(v)|| = 1.
$$
$$
|v|-1<|\beta_{3k+5}(v)| = |-v- \beta_{3k+4}|< |v| +1.
$$
For Property 2. we have
$$
\beta_{3k+4}(v) =
\frac{e_{k+2}(v)d_{3k+2}(v)d_{3k+3}(v)}{d_{k+2}(v)e_{3k+2}(v)e_{3k+3}(v)}
$$
Therefore the leading coefficient of both $e_{3k+4}(v)$ and
$d_{3k+4}(v)$ is $\pm 1$. Finally,
$$
\beta_{3k+5}(v) = -v-\beta_{3k+4}(v) =
\frac{-vd_{3k+4}(v)-e_{3k+4}(v)}{d_{3k+4}}.
$$
Since, as we have shown, the degree of $e_{3k+4}(v)$ is less than
that of $d_{3k+4}(v)$, we have that the leading coefficient of
$e_{3k+5}(v)$ comes from $-vd_{3k+4}(v)$ and therefore it equals
$\pm 1$. \endproof

By Theorem~\ref{th7}, $g_{(0,v)}(x)$ is well approximable if and
only if $v$ is a root of at least one equation $\beta_{n}(v)=0$. By
Lemma~\ref{lem9} leading coefficients of each $e_n(v)$ are plus or
minus one. Therefore all rational roots of $\beta_n(v)=0$ must also
be integers.

Assume now that $v\in\ZZ$. If $v=0$ then we obviously have
$g_{(0,0)}(x) = x^{-1}$ which is a rational function. If $
v\not\equiv 0\pmod 3$ then we use Theorem~BHWY1 for $\tilde{f}(x) =
xg_{(0,v)}(x^{-1})$. We have $C(x) = 1+vx^2, D(x)=1$ and the
functional equation~\eqref{eq_g} for $\tilde{f}(x)$ modulo 3 is
$$
(1+vx^2)\tilde(x)^2 = 1
$$
As $v\neq 0$ over $\FF_3$, we get that $\tilde{f}(x)$ is not a
rational function, therefore its continued fraction is ultimately
periodic. Going back to $\QQ$, this means that $\tilde{f}(x)$ is
badly approximable and so is $g_{(0,v)}(x)$.

Finally consider the remaining case that $v\in\ZZ, v\neq 0$ and
$v\equiv 0\pmod 3$. In this case $|v|\ge 3$ and we can use
property~3 from Lemma~\ref{lem9}. It shows that $\beta_{3k+2}(v)\neq
0$. Finally, recurrent formulae~\eqref{recur_d3} confirm that $v$ is
not a root of the remaining terms $\beta_{3k+1}(v)$ and
$\beta_{3k+3}(v)$. Application of theorem~\ref{th7} finishes the
proof.
\endproof

\noindent\textsc{Proof of Theorem~\ref{th6}}. Without loss of
generality we can assume that $u>0$. Indeed, replacing $u$ by $-u$
does not change any of the conditions~(C1),~(C2) and the property of
$g_\vu(x)$ being badly approximable is invariant under the change of
sign of $u$.

We will prove by induction that for each integer $k\ge 0$ the
following is satisfied:
\begin{equation}\label{eq_th6}
\begin{array}{l}
2u\le \alpha_{3k+2}\le 3u;\quad -2u\le \alpha_{3k+3}\le -u;\\[1ex]
|\beta_{3k+1}|\le 1;\quad u^2-v-1\le \beta_{3k+2}\le u^2-v+1;\quad
v+2u^2\le \beta_{3k+3}\le v+6u^2.
\end{array}
\end{equation}

For the base of induction we check the initial
formulae~\eqref{init_d3}: $2u\le \alpha_{2} \le 3u$ is equivalent to
$$
2(v-u^2)\le 2v-1-u^2\le 3(v-u^2).
$$
These two inequalities are in turn equivalent to $1\le u^2$ and
$v\ge 2u^2-1$ which follow from~(C1) and~(C2). The inequalities
$-2u\le \alpha_3\le -u$ follow from the fact that $\alpha_3 =
u-\alpha_2$. We obviously have $|\beta_1| = 1\le 1$ and $\beta_2 =
u^2-v \in [u^2-v-1,u^2-v+1]$. Finally we check $v+2u^2\le \beta_3\le
v+6u^2$. Since $v\ge 2u^2+8$, the enumerator $u^2+u^2+v^3-3u^2v$ of
$\beta_3$ is positive. Therefore the bounds for $\beta_3$ are
equivalent to
$$
(v+2u^2)(v-u^2)^2\le u^2+u^4+v^3-3u^2v\le (v+6u^2)(v-u^2)^2.
$$
The first inequality leads to $u^4+u^2+3u^2(u^2-1)v\ge 2u^6$ which
can easily be verified, provided that $v\ge 2u^2+8$. By simplifying
the second inequality we get
$$
u^2+u^4+11vu^4\le 3u^2+4u^2v^2+6u^6.
$$
Since $v\ge 3u^2-1$, we have $4u^2v^2\ge 4u^2(3u^2-1)v$ and the last
inequality follows from
$$
u^2+u^4+u^2v\le u^4v+6u^6
$$
which is obviously satisfied for all positive $v$ and $|u|\ge 1$.

Now we use recurrent formulae~\eqref{recur_d3} to check
estimates~\eqref{eq_th6} for $k+1$ assuming that they are satisfied
for all integer indices up to $k$.
$$
|\beta_{3k+4}| =
\left|\frac{\beta_{k+2}}{\beta_{3k+2}\beta_{3k+3}}\right|\le
\frac{v+6u^2}{(v-u^2-1)(v+2u^2)}\stackrel{(C2)}\le
\frac{v+6u^2}{(u^2+7)(v+2u^2)}.
$$
The right hand side is obviously less than one. Then the
inequalities for $\beta_{3k+5}$ follow automatically from
$\beta_{3k+5} = u^2-v-\beta_{3k+4}$ and $-1\le \beta_{3k+4}\le 1$.
Note that under conditions~(C1),~(C2), $\beta_{3k+5}$ is negative.

We have
$$
\alpha_{3k+6} =
\frac{\alpha_{k+2}+uv-\alpha_{3k+2}\beta_{3k+4}}{\beta_{3k+5}}.
$$
Since $\beta_{3k+5}<0$ inequality $-2u\le \alpha_{3k+6}$ follows
from
$$
\frac{\alpha_{k+2}+uv-\alpha_{3k+2}\beta_{3k+4}}{v-u^2-1}\le 2u
\quad\Leftarrow\quad 3u +uv +3u \le 2u(v-u^2-1)
$$
which is equivalent to $v\ge 2u^2+8$. Another inequality
$\alpha_{3k+6}\le -u$ follows from
$$
\frac{\alpha_{k+2}+uv-\alpha_{3k+2}\beta_{3k+4}}{v-u^2+1}\ge u
\quad\Leftarrow\quad -2u +uv -3u \ge u(v-u^2+1)
$$
which is equivalent to $u^2\ge 6$.

The inequalities $2u\le \alpha_{3k+5}\le 3u$ follow from the formula
$\alpha_{3k+5} = u-\alpha_{3k+6}$. Finally,
$$
v+2u^2\le v-\alpha_{3k+5}\alpha_{3k+6} \le v +6u^2
$$ which implies the inequalities for $\beta_{3k+6}$. The
claim~\eqref{eq_th6} is verified.

Inequalities~\eqref{eq_th6} suggest that $\beta_{3k+2}$ and
$\beta_{3k+3}$ can not be zeroes. The value $\beta_{3k+1}$ also can
not be zero because from~\eqref{recur_d3} it is a product of
non-zero terms. Thus the values of $\beta$ can never reach zero and
$g_\vu(x)$ is badly approximable. \endproof

\end{document}